\documentclass[12pt]{article}
\usepackage{enumerate,eepic,epic,amssymb,amsfonts,times,color,epsfig,
fancybox,graphicx,amsmath,pifont,pxfonts,txfonts,bbm}
\makeatletter
\ifcase\@ptsize
 \font\teneufm=eufm10
 \font\seveneufm=eufm7
 \font\fiveeufm=eufm5
 \font\teneusm=eusm10
 \font\seveneusm=eusm7
 \font\fiveeusm=eusm5
\or
 \font\teneufm=eufm10 scaled \magstephalf
 \font\seveneufm=eufm7
 \font\fiveeufm=eufm5
 \font\teneusm=eusm10 scaled \magstephalf
 \font\seveneusm=eusm7
 \font\fiveeusm=eusm5
\or
 \font\teneufm=eufm10 scaled \magstep1
 \font\seveneufm=eufm7
 \font\fiveeufm=eufm5
 \font\teneusm=eusm10 scaled \magstep1
 \font\seveneusm=eusm7
 \font\fiveeusm=eusm5
\fi

\newfam\eufmfam
\newfam\eusmfam
\textfont\eufmfam=\teneufm  \scriptfont\eufmfam=\seveneufm
  \scriptscriptfont\eufmfam=\fiveeufm
\textfont\eusmfam=\teneusm  \scriptfont\eusmfam=\seveneusm
  \scriptscriptfont\eusmfam=\fiveeusm

\def\frak{\ifmmode\let\next\frak@\else
 \def\next{\errmessage{Use \string\frak\space only in math mode}}\fi\next}
\def\frak@#1{{\frak@@{#1}}}
\def\frak@@#1{\fam\eufmfam#1}

\def\sh{\ifmmode\let\next\sh@\else
 \def\next{\errmessage{Use \string\sh\space only in math mode}}\fi\next}
\def\sh@#1{{\sh@@{#1}}}
\def\sh@@#1{\fam\eusmfam#1}

\ifcase\@ptsize
 \font\tenmsa=msam10
 \font\sevenmsa=msam7
 \font\fivemsa=msam5
 \font\tenmsb=msbm10
 \font\sevenmsb=msbm7
 \font\fivemsb=msbm5
\or
 \font\tenmsa=msam10 scaled \magstephalf
 \font\sevenmsa=msam7
 \font\fivemsa=msam5
 \font\tenmsb=msbm10 scaled \magstephalf
 \font\sevenmsb=msbm7
 \font\fivemsb=msbm5
\or
 \font\tenmsa=msam10 scaled \magstep1
 \font\sevenmsa=msam7
 \font\fivemsa=msam5
 \font\tenmsb=msbm10 scaled \magstep1
 \font\sevenmsb=msbm7
 \font\fivemsb=msbm5
\fi

\newfam\msafam
\newfam\msbfam
\textfont\msafam=\tenmsa  \scriptfont\msafam=\sevenmsa
  \scriptscriptfont\msafam=\fivemsa
\textfont\msbfam=\tenmsb  \scriptfont\msbfam=\sevenmsb
  \scriptscriptfont\msbfam=\fivemsb

\def\Bbb{\ifmmode\let\next\Bbb@\else
 \def\next{\errmessage{Use \string\Bbb\space only in math mode}}\fi\next}
\def\Bbb@#1{{\Bbb@@{#1}}}
\def\Bbb@@#1{\fam\msbfam#1}
\def\hexnumber@#1{\ifnum#1<10 \number#1\else
 \ifnum#1=10 A\else\ifnum#1=11 B\else\ifnum#1=12 C\else
 \ifnum#1=13 D\else\ifnum#1=14 E\else\ifnum#1=15 F\fi\fi\fi\fi\fi\fi\fi}
\def\msa@{\hexnumber@\msafam}
\def\msb@{\hexnumber@\msbfam}
\mathchardef\square="0\msa@03

\makeatother
\newcommand{\beq}{\begin{equation}}
\newcommand{\eeq}{\end{equation}}
\newcommand{\ba}{\begin{array}}
\newcommand{\ea}{\end{array}}
\newcommand{\bea}{\begin{eqnarray}}
\newcommand{\eea}{\end{eqnarray}}
\newcommand{\bean}{\begin{eqnarray*}}
\newcommand{\eean}{\end{eqnarray*}}

\newtheorem{theorem}{Theorem}[section]
\newtheorem{prop}[theorem]{Proposition}
\newtheorem{lem}[theorem]{Lemma}
\newtheorem{defi}[theorem]{Definition}
\newtheorem{cor}[theorem]{Corollary}
\newtheorem{remark}[theorem]{Remark}

\newtheorem{proof}{Proof.}
\newenvironment{pf}{\begin{proof} \rm}{\hfill$\square$ \end{proof}}

\makeatletter
\@addtoreset{equation}{section}

\makeatother

\newcommand{\rref}[1]{(\ref{#1})} 

\def\endpf{\begin{flushright}$\square$\end{flushright}}

\begin{document}
\begin{titlepage}
\begin{center}
{\huge 
Irreducible   $SL_{n+1}$--Representations   remain  Indecomposable   restricted to some Abelian Subalgebras}
\end{center}
\vspace{0.8truecm}
\begin{center}
{\large
Paolo Casati
\vskip  0.8truecm
Dipartimento di Matematica e applicazioni\\
 Universit\`a di Milano-Bicocca\\
Via degli Arcimboldi 8, I-20125 Milano, Italy\\}
\end{center}
E--mail:  casati@matapp.unimib.it\\
\vspace{0.2truecm}
\vspace{0.2truecm}
\abstract{\noindent
In this paper we show that any irreducible finite dimensional representation of $SL_{n+1}$
remains indecomposable if restricted to  $n$--dimensional abelian subalgebras 
spanned by simple root vectors.}
\vskip 1truecm\noindent{{\bf Mathematics Subjects Classification:} Primary  22E47 Secondary 17B10}\vskip 2truecm\noindent
\end{titlepage} 
\section{Introduction}
Surely the finite dimensional irreducible representations of complex simple Lie algebras are one of the most fascinating and studied subjects in the theory of representations.  Their beautiful and complicate structure  still presents unknown aspects worth to be studied (see \cite{B} and \cite{CT} for recent examples). This paper concerns with one of these, namely the restriction of such representations to  some subalgebras. More precisely we shall show that any finite dimensional irreducible representation of a   complex simple Lie algebra of type $A$ remains indecomposable if restricted to some  abelian subalgebras of maximal dimension (Theorem \ref{mthind}). Such abelian subalgebra $\mathfrak{a}$ can be constructed as follows.  Let   
$\mathfrak{g}$ be the complex simple Lie algebra  $A_n$, $\mathfrak{h}\subset\mathfrak{g}$  its Cartan subalgebra and $\Delta=\Delta(\mathfrak{g},\mathfrak{h})$ the corresponding set of roots. Further for any $\alpha \in \Delta$ let $X_\alpha$ be a basis of $\mathfrak{g}_\alpha=\{X\in \mathfrak{g}\vert \ \left[ H,X\right]=\alpha(H)X\ \forall H\in  \mathfrak{h}\}$,
$\Pi=\{\alpha_1,\dots,\alpha_n\}$ a set of simple roots in $\Delta$  and  set 
$Y_{\alpha_i}=X_{-\alpha_i}$, then $\mathfrak{a}$ is  the abelian subalgebra of  $\mathfrak{g}$ spanned by the vectors $\{Y_{\alpha_{2i+1}}\}$ ($i=0,\dots, \left[\frac{n}{2}\right]$) and 
$\{X_{\alpha_{2j}}\}$ ($j=1,\dots, \left[\frac{n}{2}\right]$), where $\left[x\right]$ denotes the integer part of $x$. \par
Theorem \ref{mthind} is  almost trivial for the Lie algebra $A_1$, while  for the Lie algebra $A_2$ was proved by Douglas and Premat in \cite{DP}, and for the remaining simple Lie algebras of rank two $B_2$ and $G_2$ by Premat in \cite{P}. 
These two papers have played an inspiring role in our work. 
As far as we know Theorem \ref{mthind} is  still unknown for $A_n$ with $n\geq 3$. \par
The paper is organized as follows. In section 2 we recall some known facts about the simple Lie algebras of type $A$ and their finite dimensional modules, and describe the maximal dimensional abelian Lie algebra $\mathfrak{a}$. This section also devoted to  present  basis of the  finite dimensional irreducible $A_n$--modules found by Littelmann in \cite{L97}. In section 3 we  find a minimal set of generators for the  restriction to the abelian subalgebra $\mathfrak{a}_n$ of  the finite dimensional representations of the Lie algebra  $A_n$, and  prove the main result of this paper: the indecomposableness  of such restricted  representations.\par
The author wishes to tank Alejandra Premat for sending her preprint \cite{P}, which plays a crucial role in the present work, and Veronica Magenes for discussions about the case concerning the Lie algebra $\mathfrak{sl}(4,\mathbbm{C})$.
\section{Irreducible finite dimensional $\mathfrak{sl}(n+1,\mathbbm{C})$--modules}
In this section we recall some basic facts on $\mathfrak{sl}(n+1,\mathbbm{C})$ and its  irreducible finite dimensional representations, and describe    the   basis  of such representations   first constructed  by Littelmann in \cite{L97}. Good references  on the structure and representation theory of the  complex simple
Lie groups and Lie algebras are,  for instance, the books \cite{H,Kn}.\par   Let $\mathfrak{g}=\mathfrak{sl}(n+1,\mathbbm{C})$ be the simple Lie algebra of all 
$(n+1)\times (n+1)$ complex matrices of  zero trace,  let $\mathfrak{h}$ be its Cartan subalgebra given by all diagonal matrices in $\mathfrak{sl}(n+1,\mathbbm{C})$,  $\mathfrak{h}^*$ its complex dual, and 
$\Delta=\Delta( \mathfrak{sl}(n+1,\mathbbm{C}),\mathfrak{h})\subset  \mathfrak{h}^*$ the corresponding set of roots.   Let $\mathfrak{g}=\mathfrak{n}^+\oplus\mathfrak{h}\oplus\mathfrak{n}^-$ be its decomposition   into the direct sum of strictly upper triangular, diagonal, and strictly lower triangular matrices, and $\Delta=\Delta^+\cup -\Delta^+$ the  decomposition of the set of root    such that 
$$
\mathfrak{n}^+= \sum_{\beta\in \Delta^+} \mathfrak{g}_\beta,\quad \mathfrak{n}^-= \sum_{\beta\in -\Delta^+} \mathfrak{g}_\beta
$$
where  $\mathfrak{g}_\beta=\{X\in \mathfrak{g}\vert\ \left[H,X\right]=\beta(H)X \ \forall H\in \mathfrak{h}\}$.
We denote by $\Pi=\{\alpha_1,\dots,\alpha_n\}$ the corresponding set of simple roots and  accordingly  we fix  a Chevalley basis of $\mathfrak{g}$: $X_\beta\in \mathfrak{g}_\beta$  and $Y_{\beta}\in \mathfrak{g}_{-\beta}$ for $\beta\in \Delta^+$, and $H_{\alpha}\in  \mathfrak{h}$ for $\alpha$ simple, in such a way that $\left[X_\alpha,Y_\alpha\right]=H_\alpha$. The Weyl group of $\mathfrak{sl}(n+1,\mathbbm{C})$ is denoted by $W$, as subgroup of $GL(\mathfrak{h}^*)$ it is generated by the hyperplane reflections 
$s_\alpha:\lambda\mapsto \lambda -\lambda(H_\alpha)\alpha$ for any $\lambda\in \mathfrak{h}^*$ and $\alpha\in \Delta$.\par
Denote by $\mathcal{U}(\mathfrak{g})$, $\mathcal{U}(\mathfrak{n}^+)$, $\mathcal{U}(\mathfrak{n}^-)$ the universal enveloping algebras of 
$\mathfrak{g}$, $\mathfrak{n}^+$, $\mathfrak{n}^-$ respectively. (More in general $\mathcal{U}(\mathfrak{a})$ will denote the universal enveloping algebra of a given subalgebra $\mathfrak{a}$  of $\mathfrak{g}$.) Following Littelmann \cite{L97} we use the following abbreviations:
$$
Y^{(k)}_\beta:=\frac{Y^k_\beta}{k!}\quad X^{(k)}_\beta:=\frac{X^k_\beta}{k!}\quad
\binom{H_\alpha}{k}:=\frac{H_\alpha(H_\alpha-1)\cdots (H_\alpha-k+1)}{k!}.
$$
Fix an ordering $\{\gamma_1 ,\dots, \gamma_N\}$  of the positive roots ($N=n(n+1)/2$). For $(\mathbf{n})\in \mathbbm{N}^N$  we set:
$$
X^{(\mathbf{n})}:=X^{(n_1)}_{\gamma_1}\cdots X^{(n_n)}_{\gamma_N},\quad Y^{(\mathbf{n})}:=Y^{(n_1)}_{\gamma_1}\cdots Y^{(n_n)}_{\gamma_N}.
$$
Fix an ordering  $\{\alpha_1 ,\dots, \alpha_n\}$ of the simple roots. For $(\mathbf{k})\in \mathbbm{N}^n$ we set:
$$
H^{(\mathbf{k})}:=\binom{H_{\alpha_1}}{k_1}\dots \binom{H_{\alpha_n}}{k_n};
$$
(we shall sometime write  $X_{i},Y_{i}, H_i$ respectively for  $X_{\alpha_i},Y_{\alpha_i},H_{\alpha_i}$, for a simple root $\alpha_i$).\par  
Recall that the monomials $Y^{(\mathbf{m})}H^{(\mathbf{k})}X^{(\mathbf{n})}$ form a Poincar\'e--Birkhoff--Witt basis of the universal enveloping algebra  $\mathcal{U}(\mathfrak{g})$,
and the monomials $X^{(\mathbf{n})}$ and $Y^{(\mathbf{m})}$ form a P--B--W basis of 
$\mathcal{U}^+=\mathcal{U}(\mathfrak{n}^+)$ respectively $\mathcal{U}^-=\mathcal{U}(\mathfrak{n}^-)$. \par  An element of $\mathfrak{h}^*$ is called a weight. The set $P=\{\lambda\in  \mathfrak{h}^*\vert \ \lambda(h_\alpha)\in \Bbb Z,\ \forall \alpha \in \Delta\}$ is said the set of  integral
weights of $\mathfrak{g}$. A weight $\lambda$  of $P$ is said dominant if 
$\lambda(H_\alpha)\geq 0$ for any simple root $\alpha$. 
The complex finite dimensional irreducible representations of $\mathfrak{sl}(n+1,\mathbbm{C})$ are parameterized by the dominant integral weights. We denote by   $V(\lambda)$ the   finite dimensional irreducible $\mathfrak{sl}(n+1,\mathbbm{C})$--module corresponding to the    integral dominant weight $\lambda$. A element $\mu$ of $\mathfrak{h}^*$ is said a weight of an  irreducible finite dimensional module $V(\lambda)$ if the weight space $V_\mu=\{v\in V(\lambda) \vert \ Hv=\mu(H)v\ \forall H\in \mathfrak{h}\}$ is different from zero. Denote  by $P(\lambda)$ the set of all weights of $V(\lambda)$ then $V(\lambda)$ may be decomposed as the direct sum
 of its weight spaces:
\beq\label{weigdec}
V(\lambda)=\bigoplus_{\mu\in P(\lambda)}V_\mu.
\eeq
Let $\Pi=\Pi_Y\cup\Pi_X$ a decomposition of the set of simple roots $\Pi$ such that the $n$--dimensional subalgebra spanned by the elements  $\{X_\alpha,Y_\beta\}_{\alpha\in \Pi_X,\beta \in \Pi_Y}$ is an abelian subalgebra. If $\alpha_1=\varepsilon_1-\varepsilon_2,\dots, \alpha_n=\varepsilon_n-\varepsilon_{n+1}$ 
is  the usual ordering of the simple roots of $\mathfrak{sl}(n+1,\mathbbm{C})$, where
$\varepsilon_i: \mathfrak{h}\to  \mathbbm{C}$ denotes the projection of a diagonal matrix onto its $i$--th entry, then 
it easy to see that this  decomposition of the set of simple roots $\Pi$  can be achieved in two  ways: either $\Pi_Y=\{\alpha_{2i+1}\}_{i=0,\dots, \left[\frac{n}{2}\right]}$ and $\Pi_X=\{\alpha_{2i}\}_{i=1,\dots, \left[\frac{n}{2}\right]}$, where $\left[x\right]$ denote the integer part of $x$ or the converse case. Since the two choices are equivalent,  let us for the sake of concreteness choose in this paper the   first one and give the  
\begin{defi}\label{an} Let $\mathfrak{a}_n$ be the abelian subalgebra of $\mathfrak{sl}(n+1,\mathbbm{C})$ spanned by the simple root--vectors $\{X_{\alpha_{2i}},Y_{\alpha_{2j+1}}\}$,  $1\leq i\leq \left[\frac{n}{2}\right]$, $0\leq j\leq \left[\frac{n}{2}\right]$.\end{defi}
\par
The aim of this paper is to show how any irreducible $\mathfrak{sl}(n+1,\mathbbm{C})$--module $V(\lambda)$  restricted to the maximal abelian subalgebra $\mathfrak{a}_n$  remains indecomposable.\par 
Further since any of such abelian algebra $\mathfrak{a}_n$ may be 
imbedded in a solvable Lie algebra endowed with  a non singular ad--invariant bilinear form \cite{MR} \cite{CP}
which is still a subalgebra of $\mathfrak{sl}(n+1,\mathbbm{C})$,
 this result provides a way to construct a fairly wide class of indecomposable (and therefore not trivial) finite dimensional modules of solvable quadratic Lie algebras  \cite{CMS}.
In order  to achieve such result   we need to consider  the basis of the irreducible $\mathfrak{sl}(n+1,\mathbbm{C})$--modules discovered by Littelmann in \cite{L97} (but see also 
\cite{L95}    \cite{L}). First we introduce the following concepts 
\begin{defi}\label{Littdef} A monomial in the $Y_{\alpha_i}$ is called semi-standard if it is of the form:
$$
Y^{(\mathbf{a})}=Y_1^{(a_1^1)}\left(Y_2^{(a_2^2)}Y_1^{(a_1^2)}\right)\bigg(\cdots\bigg)\left(Y_i^{(a_i^i)}Y_{i-1}^{(a_{i-1}^i)}\cdots
Y_1^{(a_1^i)}\right)\bigg(\cdots\bigg)\left(Y_n^{(a_n^n)}\cdots Y_2^{(a_2^n)}Y_1^{(a_1^n)}\right)
$$
where $\mathbf{a}=(a^1_1,a^2_2,a^2_1,\dots, a^n_n ,\dots, a^n_1) \in \mathbbm{N}^n$. The tuple  $\mathbf{a}$  and the monomial $Y^{(\mathbf{a})}$  are
called standard if:
{\small
$$
\mathbf{a}\in\mathcal{S}=\{(\mathbf{a})\in \mathbbm{N}^n\vert\ a^2_2\geq  a^2_1, a^3_3\geq  a^3_2\geq a^3_1,\dots,   a^i_i\geq a^i_{i-1}\geq \dots \geq a^i_1,\dots a^n_n\geq a^i_{n-1}\geq \dots\geq  a^n_1\}.
$$}
 \end{defi}
Then we can formulate the  following important result due to Littelman.
\begin{theorem}\label{Littheo} \cite{L97} For a dominant weight   $\lambda$  of $\mathfrak{g}$, let $V(\lambda)$ be  the corresponding irreducible finite dimensional $\mathfrak{g}$--module of highest weight $\lambda$
and  $u_\lambda\in V(\lambda)$ be  a highest weight vector. \\  Denote by $\lambda_i^j$ the weight of 
$$
\left( Y_i^{(a_i^j)}\cdots Y_1^{(a_1^j)}\right)\bigg(\cdots\bigg)\left(Y_n^{(a_n^n)}\cdots Y_2^{(a_2^n)}Y_1^{(a_1^n)}\right)u_\lambda
$$ 
 and set\\ $\lambda^n_0:= \lambda$,  and $\lambda^{j-1}_0:=\lambda^j_j$\ \  for $1\leq j\leq n$. 
Then the  elements of $V(\lambda)$ 
$$
Y^{(\mathbf{a})}u_\lambda=Y_1^{(a_1^1)}\left(Y_2^{(a_2^2)}Y_1^{(a_1^2)}\right)\bigg(\cdots\bigg)\left(Y_i^{(a_i^i)}Y_{i-1}^{a_{(i-1}^i)}\cdots
Y_1^{(a_1^i)}\right)\bigg(\cdots\bigg)\left(Y_n^{(a_n^n)}\cdots Y_2^{a_2^n}Y_1^{a_1^n}\right)u_\lambda
$$
with $\mathbf{a}\in \mathcal{S}$ such that 
$$
\begin{array}{llllllll}
&\lambda^n_0(H_1)\geq a^n_1 &\lambda^n_1(H_2)\geq a^n_2 &\lambda^n_2(H_3)\geq  a^n_3
 &\dots &\lambda^n_{i-1}(H_i)\geq a^n_i &\dots &\lambda^n_{n-1}(H_n)\geq a^n_n\\
&\ \ \ \ \ \  \dots\ &\ \  \ \ \ \ \ \  \dots\ &\ \ \ \ \ \ \ \ \ \dots\ & \dots\ &\ \ \ \ \ \  \ \ \ \ \ \  \dots & \dots  & \\
&\lambda^j_0(H_1)\geq a^j_1 &\ \ \ \ \ \ \ \ \dots &\ \ \ \ \ \ \ \ \ \dots 
 &\dots &\lambda^j_{j-1}(H_j)\geq a^j_j & & \\
&\ \ \ \ \ \  \dots\ &\ \  \ \ \ \ \ \  \dots\ &\ \ \ \ \ \ \ \ \ \dots\ & \dots\ &\ \ \ \ \ \  \ \ \ \ \ \   &   &\ \ \ \ \ \  \ \ \ \  \\
&\lambda^2_0(H_1)\geq a^2_1 &\lambda^2_1(H_2)\geq a^2_2 &
 & & & &\\
&\lambda^1_0(H_1)\geq a^1_1 & &
 & & & &\\
\end{array}
$$
 form  a basis $\mathfrak{L}_\lambda$  of $V(\lambda)$.
\end{theorem}
\begin{remark}\label{expfor} Let  $\lambda_i$ $i=1,\dots n$ be the elements of $\mathfrak{h}^*$ defined by the relations $\Lambda_i(\alpha_j)=\delta_{ij}$ where $\delta_{ij}$ is the usual Kronecker delta. Then if we write the dominant weight $\lambda$ in the form: $\lambda=\sum_{i=1}^mm_i\lambda_i$  (with $m_i\in \mathbbm{N}$, $i=1,\dots n$),    the conditions \rref{Littheo}  become: 
{\small
$$
\hskip -0.5truecm 
\begin{array}{llll}
&0\leq a^i_1\leq m_1-2\sum_{j=i+1}^{n}a^j_1+ \sum_{j=i+1}^{n}a^j_2&\quad  i=1,\dots, n &\\
&a^i_{k-1}\leq a^i_k\leq m_k-2\sum_{j=i+1}^{n}a^j_k+ \sum_{j=i}^{n}a^j_{k-1}+
\sum_{j=i+1}^{n}a^i_{k+1}&\quad i=1,\dots, n-k+1 &   2\leq k\leq n-1 \\
&a^n_{n-1}\leq a^n_n\leq m_n+ a^n_{n-1}. &  &\\
\end{array}
$$}
\end{remark}
Finally observe that in general the complex simple Lie algebras $\mathfrak{sl}(n+1,\mathbbm{C})$ does not have any subalgebra of dimension strictly less then $n$ with the same property of the subalgebra $\mathfrak{a}_n$ discussed above.\par\noindent
Let us indeed consider the first non trivial case, namely  the Lie algebra $\mathfrak{sl}(3,\mathbbm{C})$. In this case it is easy to show that there is not an one dimensional subalgebra such that the restriction on it of any irreducible finite dimensional representation of $\mathfrak{sl}(3,\mathbbm{C})$ remains indecomposable. Let $X$ be indeed a basis for such latter algebra. Then 
$X$ must act as a single Jordan block in any irreducible finite dimensional representation of $\mathfrak{sl}(3,\mathbbm{C})$. In particular if $\pi: \mathfrak{sl}(3,\mathbbm{C})\to \mbox{End}(\mathbbm{C}^3)$ is  the irreducible representation with $V(\lambda)=V(\lambda_1)$ (so that $\mbox{dim}_\mathbbm{C}(V(\lambda)=3$) then, since the trace of $\pi(X)$ is zero, it must exist a $\xi\in \mbox{Aut}(\mathbbm{C}^3)$ such that
$$
\xi \pi(X)\xi^{-1}=\left(\begin{array}{ccc} 0 & 0 & 0\\ 1 & 0 & 0\\ 0 & 1 & 0\end{array}\right) \qquad\mbox{or eq.}\quad  \xi \pi(X)\xi^{-1}=\left(\begin{array}{ccc} 0 & 1 & 0\\ 0 & 0 & 1\\  0 & 0 & 0\end{array}\right).
$$
I.e., $X$ can be take equal to $Y_1+Y_2$ (or eq. to $X_1+X_2$).  But the restriction to these one-dimensional sub--algebras of the module $V(\lambda)=V(\lambda_1+\lambda_2)$ is not indecomposable because on it in both cases $X^5=0$ while $\mbox{dim}(V(\lambda))=8$. 
\section{  $V(\lambda)$ as indecomposable $\mathfrak{a}_n$--module} 
Let us fix a  dominant integral weight $\lambda$ of $\mathfrak{sl}(n+1,\mathbbm{C})$. We shall show in this section that  the $\mathfrak{sl}(n+1,\mathbbm{C})$--module  $V(\lambda)$ viewed as  $\mathfrak{a}_n$--modules is indecomposable.\par We first need to find a 
(minimal) set of generators for the $\mathfrak{a}_n$--modules  $V(\lambda)$.
\begin{defi}\label{gendef}Let $V$ be a   $\mathfrak{a}_n$--modules, a subset of elements  $\{v_1,\dots v_m\}$ in $V$ is said  to be a set of generators of $V$ if $V =\mathcal{U}(\mathfrak{a}_n) \{v_1,\dots v_m\}$. The set is  called a minimal set of generators if fewer than $m$  vectors will not generate $V$. In the case of the $\mathfrak{a}_n$--modules  $V(\lambda)$ a set of generators $\mathfrak{W}$  is a set of homogeneous generators  if any element in $\mathfrak{W}$ is a $\mathfrak{sl}(n+1,\mathbbm{C})$--weight vector. \end{defi} 
\begin{theorem}\label{genset} Let $\mathcal{G}_\lambda$ be  the subset of $\mathcal{L}_\lambda=\{\mathbf{a}\in \mathcal{S}\vert \ Y^{(\mathbf{a})}u_\lambda\in \mathfrak{L}_\lambda\}$  given by: 
$$
\hskip -0.5truecm 
\mathcal{G}_\lambda =\left\{\mathbf{g}\in \mathcal{L}_\lambda\hskip 0.3truecm \vline \begin{array}{lll}  &a^{2j}_{2j}=\lambda^{2j}_{2j-1}(H_{2j}) &\ j=1,\dots\left[\frac{n}{2}\right] \\
&a^{2j+1}_{2j}\neq 0\Rightarrow a^{2j}_{2j-1}\neq 0 &\  j=1,\dots\left[\frac{n-1}{2}\right]\\
&a^{2j+1}_1=0&\  j= 0,\dots\left[\frac{n-1}{2}\right]\\
&\lambda_{2i-1}^{2j+1}(H_{2i})\neq 0 \mbox{ and } 
\lambda_0^{2j+1}(H_1)=0,\ \lambda_{2r-1}^{2j+1}(H_{2r})=a^{2j+1}_{2r}, &\  1\leq r<i \\
&\Rightarrow 
 a^{2j}_{2i-1}=\lambda_{2i-2}^{2j}(H_{2i-1})  \quad  i=1,\dots 2j-1
  &\
j=1,\dots, \left[\frac{n+1}{2}\right]
\end{array}\right\} 
$$
then  the corresponding subset $\mathfrak{G}_\lambda =\{Y^{(\mathbf{a})}u_\lambda\vert\ \  \mathbf{a}\in \mathcal{G}_\lambda\}$ of the Littelmann basis $\mathfrak{L}_\lambda$ 
 is a set of homogeneous $\mathfrak{a}_n$--generators of $V(\lambda)$.
\end{theorem}
{\bf Proof} Set $\mathcal{L}^0_\lambda=\{\mathbf{a}\in \mathcal{L}_\lambda\vert \ \  a^1_1=0\}$. Since $X_1$ belongs to  $\mathfrak{a}_n$, we have of course  only to prove that acting with  $\mathfrak{a}_n$ we may construct the subset  $\mathfrak{L}^0_\lambda=\{Y^{(\mathbf{a})}u_\lambda\vert  \mathbf{a}\in \mathcal{L}^0_\lambda \}$  of $\mathfrak{L}_\lambda$.  We divide the proof  in four  steps.\par\noindent
1.  First,  if we define 
$$
\hskip -1truecm
\mathcal{L}^1_\lambda=\left\{\mathbf{a}\in \mathcal{L}^0_\lambda\vline \begin{array}{lll}  &a_1^{2h+1}=0,\  a^{2h+1}_{2h}\neq 0\Rightarrow a^{2h}_{2h-1}\neq 0 &\  h= 1,\dots\left[\frac{n-1}{2}\right]\\
&\lambda_{2i-1}^{2j+1}(H_{2i})\neq 0 \mbox{ and } 
\lambda_0^{2j+1}(H_1)=0,\ \lambda_{2r-1}^{2j+1}(H_{2r})=a^{2j+1}_{2r}, &\  1\leq r<i \\
&\Rightarrow 
 a^{2j}_{2i-1}=\lambda_{2i-2}^{2j}(H_{2i-1})  \quad  i=1,\dots 2j-1
  &\
j=1,\dots, \left[\frac{n+1}{2}\right]\\
\end{array}
\right\}
$$
and set  $\mathfrak{L}^1_\lambda=\{Y^{(\mathbf{a})}\vert\   \mathbf{a}\in \mathcal{L}^1_\lambda\}$, then  
$\mathfrak{L}^1_\lambda\subset\mathcal{U}(\mathfrak{a}_n)(\mathfrak{G}_\lambda)$.\par  
Let us  consider indeed for any $1\leq j\leq \left[\frac{n}{2}\right]$ the  subsets  $\mathfrak{L}^1_{2j}$ of  $\mathfrak{L}^1_\lambda$ given by:
{\small
$$
\hskip -0.5truecm \mathfrak{L}^1_{2j}=\left\{Y^{(\mathbf{a})}u_\lambda\in \mathfrak{L}^1_\lambda \vert\ \exists\ Y^{(\mathbf{a}_{\mathcal{G}})}u_\lambda\in \mathfrak{G} \
\vline\begin{array}{llll}
  (a_{\mathcal{G}})^h_{2l+1}=a_{2l+1}^{h},& &   l=0,\dots, \left[\frac{n-1}{2}\right],  &2l+1\leq  h\leq n \\
(a_{\mathcal{G}})^h_{2k}=a_{2k}^{h},& &
 k=j,\dots, \left[\frac{n}{2}\right] &2l\leq  h\leq n, 
\end{array}\right\}
$$}
and the corresponding filtration of $\mathfrak{L}^1_\lambda$: 
$$
\mathfrak{G}_\lambda=\mathfrak{L}^1_2\subset\dots\subset\mathfrak{L}^1_{2j}\dots\subset
\mathfrak{L}^1_{2\left[\frac{n}{2}\right]}\subset\mathfrak{L}^1_{2\left[\frac{n}{2}\right]+2} =\mathfrak{L}^1_\lambda.
$$
 Obviously it suffices to show that $\mathfrak{L}^1_{2j}\subset \mathcal{U}(\mathfrak{a}_n)(\mathfrak{L}^1_{2j-2})$ for any $1\leq j\leq \left[\frac{n}{2}\right]$. We shall do 
   it  (for a fixed index $j$) by induction over the partial ordering ``$\leq_{j}$'' of $\mathcal{L}^0_\lambda$ (and of  $\mathfrak{L}^0_\lambda$  as well) given by the relations
$$
\mathbf{a}\leq_{j}\mathbf{b} \Leftrightarrow  a^i_{2j}-a^i_{2j+1} \leq 
 b^i_{2j}-b^i_{2j+1} \quad i=2j+1,\dots, n.
$$
With respect to this ordering the minimal elements in $\mathfrak{L}^1_{2j}$ are those  $Y^{(\mathbf{a})}u_\lambda$ with $a^i_{2j}-a^i_{2j+1}=-\lambda^i_{2j}(H_{2j+1})$, $2j+1\leq i\leq n$. For any of this element there exists a positive integer number $k$ (namely  $k=a^{2j}_{2j}-\lambda_{2j+1}^{2j}(H_{2j})$) such that the element $Y^{(\mathbf{g(a)})}u_\lambda$  with ${g(a)}^{2j}_{2j}=a^{2j}_{2j}+k$ and  ${g(a)}^{i}_{l}=a^{i}_{l}$ if $(i,l)\neq (2j,2j)$ belongs to $\mathfrak{G}_\lambda$ 
 and ( using the  $\mathfrak{sl}(2,\mathbbm{C})$--representation theory and \cite{L97}) 
$$
X_{2j}^k\cdot Y^{(\mathbf{g(a)})}u_\lambda=c^{2j}_k(\mathbf{a}) Y^{(\mathbf{a})}u_\lambda=
\prod_{i=1}^{k-1}\left(\sum_{h=1}^{n-2j}a^{2j+h}_{2j+1}-2\sum_{h=1}^{n-2j}a^{2j+h}_{2j}+\sum_{h=0}^{n-2j}a^{2j+h}_{2j-1}-a^{2j}_{2j}+k+i\right)
 Y^{(\mathbf{a})}u_\lambda$$
hence, since in our Hypothesis the coefficients $c^{2j}_k(\mathbf{a})$ are always different from zero,  any minimal element of $\mathfrak{L}^1_{2j}$
belongs to   $\mathcal{U}(\mathfrak{a}_n)(\mathfrak{L}^1_{2j-2})$.  Suppose now by induction hypothesis that we have constructed any elements  $Y^{(\mathbf{b})}u_\lambda\in \mathfrak{L}^1_{2j}$  for any $\mathbf{b}<_{j}\mathbf{a}$. Since there exists a tuple  $\mathbf{a}_{\mathcal{L}^1_{2j-2}}$ such that $(a_{\mathcal{L}_{2j-2}^1})^i_{l}=a^i_{l}$, if $(i,l)\neq (2j,2j)$  
 and $Y^{(\mathbf{a}_{\mathcal{L}^1_{2j-2}})}u_\lambda\in \mathfrak{L}^1_{2j-2}$,
we have 
 $$
X_{2j}^k\cdot Y^{(\mathbf{a}_{\mathcal{L}^1_{2j-2}})}u_\lambda=c^{2j}_k(\mathbf{a}_{\mathcal{L}^1_{2j-2}})
Y^{(\mathbf{a})}u_\lambda+\sum_{\mathbf{b}<_{j}\mathbf{a}}c_\mathbf{b}Y^{(\mathbf{b})}u_\lambda
$$
which shows (being again $c^{2j}_k(\mathbf{a}_{\mathcal{L}^1_{2j-2}})\neq 0$)
that also $Y^{(\mathbf{a})}u_\lambda$ belongs to $\mathcal{U}(\mathfrak{a}_n)(\mathfrak{L}^1_{2j-2})$.\par\noindent  
2. Define  now  
$$
\hskip -1truecm
\mathcal{L}^2_\lambda=\left\{\mathbf{a}\in \mathcal{L}^0_\lambda\vline \begin{array}{lll}  &a_1^{2h+1}=0, &\  h= 1,\dots\left[\frac{n-1}{2}\right]\\
&\lambda_{2i-1}^{2j+1}(H_{2i})\neq 0 \mbox{ and } 
\lambda_0^{2j+1}(H_1)=0,\ \lambda_{2r-1}^{2j+1}(H_{2r})=a^{2j+1}_{2r}, &\ 1\leq r<i \\
&\Rightarrow 
 a^{2j}_{2i-1}=\lambda_{2i-2}^{2j}(H_{2i-1})  \quad  i=1,\dots 2j-1
  &\
j=1,\dots, \left[\frac{n+1}{2}\right]
\end{array}\right\} 
$$
and  set   $\mathfrak{L}^2_\lambda=\left\{Y^{(\mathbf{a})}u_\lambda\vert\ \mathbf{a}\in \mathcal{L}^2_\lambda \right\})$  
then $\mathfrak{L}^2_\lambda\subset\mathcal{U}(\mathfrak{a}_n)(\mathfrak{L}^1_\lambda)$. \par
For any $1\leq s\leq \left[\frac{n-1}{2}\right]$ let $\mathfrak{L}^2_s$ be  the set 
$$
\mathfrak{L}^2_s=\left\{Y^{(\mathbf{a})}u_\lambda\in \mathfrak{L}^2_\lambda\vert\ a^{2k+1}_{2k}\neq 0\Rightarrow a^{2k}_{2k-1}\neq 0\quad k=s,\dots \left[\frac{n-1}{2}\right]\right\}
$$ 
and consider the corresponding filtration of $\mathfrak{L}^2_\lambda$: 
$$
\mathfrak{L}^1_\lambda=\mathfrak{L}^2_1\subset \dots\subset \mathfrak{L}^2_s\subset \dots\subset \mathfrak{L}^2_{\left[\frac{n-1}{2}\right]}\subset \mathfrak{L}^2_{\left[\frac{n-1}{2}\right]+1}=\mathfrak{L}^2_\lambda.
$$
 Again it will suffice to show that $\mathfrak{L}^2_s\subset  \mathcal{U}(\mathfrak{a}_n)(\mathfrak{L}^1_{s-1})$ for any $2\leq s\leq \left[\frac{n-1}{2}\right]+1$.
We shall still do it  by induction. Indeed consider first an element in $\mathfrak{L}^2_s$ of the form $Y^{(\mathbf{a})}u_\lambda
=\left(Y_2^{(a^2_2)}\cdots\right)\bigg(\cdots\bigg)\left(Y_{2s}^{(a^{2s}_{2s})}
Y_{2s-2}^{(a^{2s}_{2s-2})}\cdots\right)\left(Y_{2s+1}^{(a^{2s+1}_{2s+1})}Y_{2s}\cdots\right)
\bigg(\cdots\bigg)u_\lambda$ then the element  $Y^{(\mathbf{b(a)})}u_\lambda$ with the tuple $\mathbf{b(a)}$ given by the relations $b(a)^{2s}_{2s}={a}^{2s}_{2s}+1$,   $b(a)^{2s+1}_{2s+1}={a}^{2s+1}_{2s+1}-1$,  $b(a)^{2s+1}_{2s}=0$ and $b(a)^{i}_{j}={a}^{i}_{j}$ otherwise, belongs to $\mathfrak{L}^2_{s-1}$. Further from the relation   \cite{L97}
$$
\hskip -0.5truecm
\begin{array}{ll}
&Y_{2s+1}\cdot Y^{(\mathbf{b(a)})}u_\lambda\\
&=p(1,{a}^{2s}_{2s}+1,{a}^{2s+1}_{2s+1},0)\bigg(\cdots\bigg)\left( Y^{(a^{2s}_{2s}+1)}_{2s}Y^{(a^{2s}_{2s-2})}_{2s-2}\cdots\right)\left(Y^{(a^{2s+1}_{2s+1})}_{2s+1}Y^{(a^{2s+1}_{2s+2})}_{2s-1}
\cdots\right)\left(Y^{(a^{2s+2}_{2s+2})}_{2s+2}
\cdots\right)\bigg(\cdots\bigg)u_\lambda\\&+p(1,{a}^{2s}_{2s}+1,{a}^{2s+1}_{2s+1},1)\bigg(\cdots\bigg)\left(Y_{2s}^{(a^{2s}_{2s})}
Y_{2s-2}^{(a^{2s}_{2s-2})}\cdots\right)\left(Y_{2s+1}^{(a^{2s+1}_{2s+1})}Y_{2s}\cdots\right)\left(Y^{(a^{2s+2}_{2s+2})}_{2s+2}
\cdots\right)\bigg(\cdots\bigg)u_\lambda,
\end{array}
$$
where 
$$
p(a,b,c,d)=\binom{a+c-b}{a-d}\quad a,b,c,d\in \mathbbm{N}\quad d\leq b
$$
and to have binomial coefficients also available for negative integers, following Littelmann  we used the definition:
$$
\binom{a}{b}=\lim_{t\to 0} \frac{\Gamma(a+1+t)}{\Gamma(b-a+1+t)\Gamma(b+1+t)};
$$
it follows   that $Y^{(\mathbf{a})}u_\lambda$ belongs to 
$\mathcal{U}(\mathfrak{a}_n)\left(\mathfrak{L}^2_{s-1}\right)$ because $p(1,{a}^{2s}_{2s}+1,{a}^{2s+1}_{2s+1},1)=1$ (but also see \cite{L97} remark 7) and both $Y^{(\mathbf{b(a)})}u_\lambda$ and \\ $\bigg(\cdots\bigg)\left(Y^{(a^{2s}_{2s})}_{2s}\cdots\right)\left(Y^{(a^{2s+1}_{2s+1})}_{2s+1}Y^{(a^{2s+1}_{2s+2})}_{2s-1}
\cdots\right)\left(Y^{(a^{2s+2}_{2s+2})}_{2s+2}\cdots\right)\bigg(\cdots\bigg)u_\lambda$ are in $\mathcal{U}(\mathfrak{a}_n)\left(\mathfrak{L}^2_{s-1}\right)$.\\ Let us now consider an element in $\mathcal{L}^2_s$ of the type \\
$\left(Y_2^{(a^2_2)}\cdots\right)\bigg(\cdots\bigg)\left(Y_{2s}^{(a^{2s}_{2s})}
Y_{2s-2}^{(a^{2s}_{2s-2})}\cdots\right)\left(Y_{2s+1}^{(a^{2s+1}_{2s+1})}Y^{k+1}_{2s}\cdots\right)
\left(Y_{2s+2}^{(a^{2s+2}_{2s+2})}\right)\bigg(\cdots\bigg)u_\lambda$, since by induction\\ Hypothesis $Y^{(\mathbf{b})}u_\lambda$, with  $b^{2s}_{2s}=a^{2s}_{2s}+1$,   $b^{2s+1}_{2s+1}=a^{2s+1}_{2s+1}-1$,  $b^{2s+1}_{2s}=k$ and ${b}^{i}_{j}={a}^{i}_{j}$ otherwise, 
belongs to $\mathcal{U}(\mathfrak{a}_n)\left(\mathcal{L}^2_{s-1}\right)$ from  
$$
\hskip -0.3truecm
\begin{array}{ll}
&Y_{2s+1}\cdot Y^{(\mathbf{b})}u_\lambda\\
&=p(1,{a}^{2s}_{2s}+1,{a}^{2s+1}_{2s+1},0)\bigg(\cdots\bigg)\left( Y^{(a^{2s}_{2s}+1)}_{2s}\cdots\right)\left(Y^{(a^{2s+1}_{2s+1})}_{2s+1}Y^{(a^{2s+1}_{2s+2})}_{2s-1}Y_{2s}^k
\cdots\right)\left(Y^{(a^{2s+2}_{2s+2})}_{2s+2}
\cdots\right)\bigg(\cdots\bigg)u_\lambda\\&+p(1,{a}^{2s}_{2s}+1,{a}^{2s+1}_{2s+1},1)\bigg(\cdots\bigg)\left(Y_{2s}^{(a^{2s}_{2s})}
Y_{2s-2}^{(a^{2s}_{2s-2})}\cdots\right)\left(Y_{2s+1}^{(a^{2s+1}_{2s+1})}Y^{k+1}_{2s}\cdots\right)\left(Y^{(a^{2s+2}_{2s+2})}_{2s+2}
\cdots\right)\bigg(\cdots\bigg)u_\lambda
\end{array}
$$
it follows that also $\left(Y_2^{(a^2_2)}\cdots\right)\bigg(\cdots\bigg)\left(Y_{2s}^{(a^{2s}_{2s})}
Y_{2s-2}^{(a^{2s}_{2s-2})}\cdots\right)\left(Y_{2s+1}^{(a^{2s+1}_{2s+1})}Y^{k+1}_{2s}\cdots\right)
\left(Y_{2s+2}^{(a^{2i+2}_{2i+2})}\right)\bigg(\cdots\bigg)u_\lambda$
  belongs to $\mathcal{U}(\mathfrak{a}_n)\left(\mathcal{L}^2_{s-1}\right)$. \par\noindent
3. Let us now define   $\mathcal{L}^3_\lambda$ as:
$$
\hskip -1truecm
\mathcal{L}^3_\lambda=\left\{\mathbf{a}\in \mathcal{L}_0\vline  \begin{array}{lll} &\lambda_{2i-1}^{2j+1}(H_{2i})\neq 0 \mbox{ and } 
\lambda_0^{2j+1}(H_1)=0,\ \lambda_{2r-1}^{2j+1}(H_{2r})=a^{2j+1}_{2r}, &\ 1\leq r<i \\
&\Rightarrow 
 a^{2j}_{2i-1}=\lambda_{2i-2}^{2j}(H_{2i-1})  \quad  i=1,\dots 2j-1
  &\
j=1,\dots, \left[\frac{n+1}{2}\right]
\end{array}\right\} 
$$ 
and set $\mathfrak{L}^3_\lambda=\{Y^{(\mathbf{a})}u_\lambda\vert\ \mathbf{a}\in \mathcal{L}^3_\lambda\}$ then  $\mathfrak{L}^3_\lambda\subset \mathcal{U}(\mathfrak{a}_n)( \mathfrak{L}^2_\lambda)$.\par\noindent
Defining   for any $1\leq s\leq  \left[\frac{n-1}{2}\right]$ the sets 
$$
\mathfrak{L}^3_s=\{Y^{(\mathbf{a})}u_\lambda\in \mathfrak{L}^3_\lambda\vert\ a^{2j+1}_1=0\quad j=s,\dots  \left[\frac{n-1}{2}\right]\}
$$
we have  the  filtration of $\mathcal{L}^3$:
$$
\mathfrak{L}^2_\lambda=\mathfrak{L}^3_1\subset \dots\subset \mathfrak{L}^3_s\subset \dots\subset \mathfrak{L}^3_{\left[\frac{n-1}{2}\right]}\subset   \mathfrak{L}^3_{\left[\frac{n-1}{2}\right]+1}=\mathfrak{L}^3_\lambda.
$$
Once again it suffices to  prove that  $\mathfrak{L}^3_s\subset\mathcal{U}(\mathfrak{a}_n)(\mathfrak{L}^3_{s-1})$ for any fixed $s$. We proceed  by induction.  Let us first consider an element 
of $\mathfrak{L}^3_s$ of the type \\
$Y^{(\mathbf{a})}u_\lambda
=\left(Y_2^{(a^2_2)}\cdots\right)\bigg(\cdots\bigg)\left(Y_{2s}^{(a^{2s}_{2s})}\cdots
Y_{1}^{(a^{2s}_{1})}\right)\left(Y_{2s+1}^{(a^{2s+1}_{2s+1})}
\dots Y_{2}^{(a^{2s+1}_{2})}Y_{1}\right)
\bigg(\cdots\bigg)u_\lambda$ then the element $Y^{(\mathbf{b(a)})}u_\lambda$
with $b(a)^{2s}_{2s-l}=a^{2s}_{2s-l}+1$, $0\leq l\leq 2s-1$ ,  $b(a)^{2s+1}_{2s-l}=a^{2s+1}_{2s-l}-1$, $-1\leq l\leq 2s-1$, $b(a)^h_l=a^h_l$ otherwise, belongs to $\mathfrak{L}^3_{s-1}$, and we have 
{\small
$$
\hskip -0.5truecm
\begin{array}{ll}
&Y_{2s+1}\cdot Y^{(\mathbf{b(a)})}u_\lambda=p(1,a^{2s}_{2s}+1, a^{2s+1}_{2s+1}-1,0)\bigg(\cdots\bigg)\left(Y_{2s}^{(a^{2s}_{2s}+1)}\cdots
\right)\left(Y_{2s+1}^{(a^{2s+1}_{2s+1})}Y_{2s}^{(a^{2s+1}_{2s}-1)}
\cdots \right)\bigg(\cdots\bigg)u_\lambda\\&+ p(1,a^{2s}_{2s}+1, a^{2s+1}_{2s+1}-1,1)\bigg(\cdots\bigg)\left(Y_{2s}^{(a^{2s}_{2s})}
Y_{2s+1}^{(a^{2s+1}_{2s+1})}Y_{2s}Y_{2s-1}^{(a^{2s}_{2s-1}+1)}\cdots\right)\left(
 Y_{2s}^{(a^{2s+1}_{2s}-1)}\cdots\right)\bigg(\cdots\bigg)u_\lambda\\
&=\dots\\
&=\sum_{k=1}^h p^s_k(\mathbf{a})\bigg(\cdots\bigg)\left(\prod_{l=0}^{k-2}Y_{2s-l}^{(a^{2s}_{2s-l})}\right)
\left(\prod_{l=k-1}^{2s-1}Y_{2s-l}^{(a^{2s}_{2s-l}+1)}\right)\left(\prod_{l=0}^{k-1}Y_{2s+1-l}^{(a^{2s+1}_{2s+1-l})}\right)
\left(\prod_{l=k}^{2s-1}Y_{2s+1-l}^{(a^{2s+1}_{2s+1-l}-1)}\right)\bigg(\cdots\bigg)u_\lambda\\&+ 
q^s_h(\mathbf{a})\bigg(\cdots\bigg)\left(\prod_{l=0}^{h-1}Y_{2s-l}^{(a^{2s}_{2s-l})}\right)\left(\prod_{l=0}^{h-1}Y_{2s+1-l}^{(a^{2s+1}_{2s+1-l})}\right)
\left(Y_{2s-h+1}Y_{2s-h}^{(a^{2s}_{2s-h}+1)}
Y_{2s-h+1}^{(a^{2s+1}_{2s-h+1}-1)}\right)\\&
\left(\prod_{l=h+1}^{2s-1}Y_{2s-l}^{(a^{2s}_{2s-l}+1)}\right)
\left(\prod_{l=h+1}^{2s-1}Y_{2s+1-l}^{(a^{2s+1}_{2s+1-l}-1)}\right)\bigg(\cdots\bigg)u_\lambda=\dots\\& 
=\sum_{k=1}^{2s} p^s_k(\mathbf{a})\left(\cdots\right)\left(\prod_{l=0}^{k-2}Y_{2s-l}^{(a^{2s}_{2s-l})}\right)
\left(\prod_{l=k-1}^{2s-1}Y_{2s-l}^{(a^{2s}_{2s-l}+1)}\right)\left(\prod_{l=0}^{k-1}Y_{2s+1-l}^{(a^{2s+1}_{2s+1-l})}\right)
\left(\prod_{l=k}^{2s-1}Y_{2s+1-l}^{(a^{2s+1}_{2s+1-l}-1)}\right)\bigg(\cdots\bigg)u_\lambda
\\&+q^s_{2s}(\mathbf{a})Y^{(\mathbf{a})}u_\lambda
\end{array}
$$}
 where 
$$
\begin{array}{ll}
p^s_k(\mathbf{a})&=p(1,a^{2s}_{2s+1-k}+1, a^{2s+1}_{2s-k+2}-1,0)\prod_{l=1}^{k-1}p(1,a^{2s}_{2s+1-l}+1, a^{2s+1}_{2s-l+2}-1,1)\\
q^s_k(\mathbf{a})&=\prod_{l=1}^{k}p(1,a^{2s}_{2s+1-l}+1, a^{2s+1}_{2s-l+2}-1,1),
\end{array}
$$
 which show that $Y^{(\mathbf{a})}u_\lambda\in \mathcal{U}(\mathfrak{a}_n)(\mathfrak{L}^3_{s-1})$ because $q^s_{2s}(\mathbf{a})\neq 0$ and all other elements belong to  $\mathcal{U}(\mathfrak{a}_n)(\mathfrak{L}^3_{s-1})$. If  we now  consider 
an element $ Y^{(\mathbf{a})}u_\lambda$ of the type \\
$Y^{(\mathbf{a})}u_\lambda
=\left(Y_2^{(a^2_2)}\cdots\right)\bigg(\cdots\bigg)\left(Y_{2s}^{(a^{2s}_{2s})}\cdots
\dots Y_{1}^{(a^{2s}_{1})}\right)\left(Y_{2s+1}^{(a^{2s+1}_{2s+1})}
\dots Y_{2}^{(a^{2s+1}_{2})}Y^k_{1}\right)
\bigg(\cdots\bigg)u_\lambda$  then by induction Hypothesis the element $Y^{(\mathbf{b(a)})}u_\lambda$
with $b(a)^{2s}_{2s-l}=a^{2s}_{2s-l}+1$, $0\leq l\leq 2s-1$ ,  $b(a)^{2s+1}_{2s-l}=a^{2s+1}_{2s-l}-1$, $-1\leq l\leq 2s-1$, $b(a)^h_l=a^h_l$ otherwise, belongs to $\mathfrak{L}^3_{s-1}$, and with the same computations done before  we have
{\small $$
\hskip -0.4truecm
\begin{array}{ll}
&Y_{2s+1}\cdot Y^{(\mathbf{b(a)})}u_\lambda\\
&=\sum_{j=1}^{2s} p^s_j(\mathbf{a})\bigg(\cdots\bigg)\left(\Pi_{l=0}^{j-2}Y_{2s-l}^{(a^{2s}_{2s-l})}\right)
\left(\prod_{l=j-1}^{2s-1}Y_{2s-l}^{(a^{2s}_{2s-l}+1)}\right)\left(\prod_{l=0}^{j-1}Y_{2s+1-l}^{(a^{2s+1}_{2s+1-l})}\right)
\left(\prod_{l=j}^{2s-1}Y_{2s+1-l}^{(a^{2s+1}_{2s+1-l}-1)}\right)\bigg(\cdots\bigg)u_\lambda
\\&+q^s_{2s}(\mathbf{a})Y^{(\mathbf{a})}u_\lambda,
\end{array}
$$}
which  implies  (by induction) $Y^{(\mathbf{a})}u_\lambda\in \mathcal{U}(\mathfrak{a}_n)(\mathfrak{L}^3_{s-1})$.
\par\noindent
4. Finally we can show that  $\mathfrak{L}^0_\lambda\subset\mathcal{U}(\mathfrak{a}_n)(\mathfrak{L}^3_\lambda)$. The computations are similar to those done in the previous step.
We define indeed for any $1\leq s\leq \left[\frac{n-1}{2}\right]$ the sets 
{\small
$$
\mathfrak{L}^0_s=\left\{Y^{(\mathbf{a})}u_\lambda\in \mathfrak{L}^0_\lambda\vline 
\begin{array}{lll}
&\lambda_{2i-1}^{2j+1}(H_{2i})\neq 0 \mbox{ and } 
\lambda_0^{2j+1}(H_1)=0,\ \lambda_{2r-1}^{2j+1}(H_{2r})=a^{2j+1}_{2r}, &\ 1\leq r<i\\ 
&\Rightarrow
 a^{2j}_{2i+1}=\lambda_{2i}^{2j}(H_{2i+1}) \quad    i=1,\dots 2j-1
&\ 
j=s,\dots, \left[\frac{n+1}{2}\right]
\end{array}\right\}
$$} 
and consider the corresponding filtration of $\mathfrak{L}^0_\lambda$
$$
\mathfrak{L}^3_\lambda=\mathfrak{L}^0_1\subset \dots\subset \mathfrak{L}^0_s\subset \dots\subset \mathfrak{L}^0_{\left[\frac{n+1}{2}\right]}\subset   \mathfrak{L}^0_{\left[\frac{n+1}{2}\right]+1}=\mathfrak{L}^0_\lambda.
$$
Again, we need only to prove (always  by induction)   that $\mathfrak{L}^0_s\subset \mathcal{U}(\mathfrak{a}_n)( \mathfrak{L}^0_{s-1})$ for any $1\leq s\leq \left[\frac{n-1}{2}\right]$.\\
 For a fixed $i$, $1\leq i\leq 2j-1$, let  $Y^{(\mathbf{a})}u_\lambda\in  \mathfrak{L}^0_s$ be of the type  
$$
Y^{(\mathbf{a})}u_\lambda=\left(Y_2^{(a^2_2)}\cdots\right)\bigg(\cdots\bigg)\left(Y_{2s}^{(a^{2s}_{2s})}
\cdots Y_1^{(a^{2s}_1)}\right)\left(Y_{2s+1}^{(a^{2s+1}_{2s+1})}\cdots Y_{2i}^{(\lambda^{2s+1}_{2i-2}(H_{2i-1})+a^{2s}_{2i-1}+1)}\cdots \right)\bigg(\cdots\bigg)u_\lambda
$$  
then the element $Y^{(\mathbf{b(a)})}u_\lambda$
with $b(a)^{2s}_{2s-l}=a^{2s}_{2s-l}+1$, $0\leq l\leq 2s-2i$, $b(a)^{2s+1}_{2s+1-l}=a^{2s+1}_{2s+1-l}-1$, $0\leq l\leq 2s-2i+1$, $b(a)^h_l=a^h_l$ otherwise, belongs to $\mathfrak{L}_{s-1}^0$ and we have with the same computations of the previous step and the results of \cite{L97}:
{\small $$
\hskip -0.5truecm
\begin{array}{ll}
&Y_{2s+1}\cdot Y^{(\mathbf{b(a)})}u_\lambda\\
&=\sum_{k=1}^{2s-2i+1} p^s_k(\mathbf{a})\bigg(\cdots\bigg)\left(\Pi_{l=0}^{k-2}Y_{2s-l}^{(a^{2s}_{2s-l})}\right)
\left(\Pi_{l=k-1}^{2s-1}Y_{2s-l}^{(a^{2s}_{2s-l}+1)}\right)\left(\Pi_{l=0}^{k-1}Y_{2s+1-l}^{(a^{2s+1}_{2s+1-l})}\right)
\left(\Pi_{l=k}^{2s-1}Y_{2s+1-l}^{(a^{2s+1}_{2s+1-l}-1)}\right)\bigg(\cdots\bigg)u_\lambda
\\&+p^s_{2s-2i+2}(\mathbf{a})Y^{(\mathbf{a})}u_\lambda\\
&+\sum_{k=s-i+2}^sp^s_{2k-1}(\mathbf{a})\bigg(\cdots\bigg)\left(\Pi_{l=0}^{2k-3}Y_{2s-l}^{(a^{2s}_{2s-l})}\right)
\left(\Pi_{l=2k-1}^{2s-1}Y_{2s-l}^{(a^{2s}_{2s-l}+1)}\right)\left(\Pi_{l=0}^{2k-2}Y_{2s+1-l}^{(a^{2s+1}_{2s+1-l})}\right)\\
&\left(\Pi_{l=2k-1}^{2s-1}Y_{2s+1-l}^{(a^{2s+1}_{2s+1-l}-1)}\right)\bigg(\cdots\bigg)u_\lambda
\end{array}
$$}
which implies that $Y^{(\mathbf{a})}u_\lambda$  belongs to $\mathcal{U}(\mathfrak{a}_n)(\mathfrak{L}_{s-1}^0)$.   
Now suppose by induction Hypothesis that we have already constructed all the elements of $\mathfrak{L}_s^0$ with $a^{2s+1}_{2i}=\lambda_{2i-2}^{2s+1}(H_{2i-1})+a^{2s}_{2i-1}+k$. Then for any element $Y^{(\mathbf{a})}u_\lambda\in \mathfrak{L}^0_s$  of the type $$
Y^{(\mathbf{a})}u_\lambda=\left(Y_2^{(a^2_2)}\cdots\right)\bigg(\cdots\bigg)\left(Y_{2s}^{(a^{2s}_{2s})}
\cdots Y_1^{(a^{2s}_1)}\right)\left(Y_{2s+1}^{(a^{2s+1}_{2s+1})}\cdots Y_{2i}^{(\lambda^{2s+1}_{2i-2}(H_{2i-1})+a^{2s-1}_{2i-1}+k+1)}\cdots \right)\bigg(\cdots\bigg)u_\lambda
$$ 
the element $Y^{(\mathbf{b(a)})}u_\lambda$
with $b(a)^{2s}_{2s-l}=a^{2s}_{2s-l}+1$, $0\leq l\leq 2s-2i$, $b(a)^{2s+1}_{2s+1-l}=a^{2s+1}_{2s+1-l}-1$, $0\leq l\leq 2s-2i+1$, $b(a)^h_l=a^h_l$ otherwise, belongs to $\mathfrak{L}_{s-1}^0$ and we 
have
{\small $$
\hskip -0.5truecm
\begin{array}{ll}
&Y_{2s+1}\cdot Y^{(\mathbf{b(a)})}u_\lambda\\
&=\sum_{k=1}^{2s-2i+1} p^s_k(\mathbf{a})\bigg(\cdots\bigg)\left(\Pi_{l=0}^{k-2}Y_{2s-l}^{(a^{2s}_{2s-l})}\right)
\left(\Pi_{l=k-1}^{2s-1}Y_{2s-l}^{(a^{2s}_{2s-l}+1)}\right)\left(\Pi_{l=0}^{k}Y_{2s+1-l}^{(a^{2s+1}_{2s+1-l})}\right)
\left(\Pi_{l=k+1}^{2s-1}Y_{2s+1-l}^{(a^{2s+1}_{2s+1-l}-1)}\right)\bigg(\cdots\bigg)u_\lambda
\\&+p^s_{2s-2i+2}(\mathbf{a})Y^{(\mathbf{a})}u_\lambda\\
&+\sum_{k=s-i+2}^sp^s_{2k-1}(\mathbf{a})\bigg(\cdots\bigg)\left(\Pi_{l=0}^{2k-3}Y_{2s-l}^{(a^{2s}_{2s-l})}\right)
\left(\Pi_{l=2k-1}^{2s-1}Y_{2s-l}^{(a^{2s}_{2s-l}+1)}\right)\left(\Pi_{l=0}^{2k-2}Y_{2s+1-l}^{(a^{2s+1}_{2s+1-l})}\right)\\
&\left(\Pi_{l=2k-1}^{2s-1}Y_{2s+1-l}^{(a^{2s+1}_{2s+1-l}-1)}\right)\bigg(\cdots\bigg)u_\lambda
\end{array}
$$}
which once again implies that $Y^{(\mathbf{b(a)})}u_\lambda\in \mathcal{U}(\mathfrak{a}_n)(\mathfrak{L}_{s-1}^0)$. This closes the proof of the Theorem
\endpf
In the first non trivial case beyond that treated by Douglas and Premat \cite{DP}    namely the restriction of irreducible finite dimensional $\mathfrak{sl}(4,\mathbbm{C})$--module $V(\lambda)$, $\lambda=n\Lambda_1+m\Lambda_2+p\Lambda_3$ to the abelian three dimensional Lie algebra spanned by the element $Y_1,X_2,Y_3$ the set of generators $\mathfrak{G}_\lambda$ is:
$$
\mathfrak{G}_{\lambda}=\left\{Y_2^{m-j+i+h}Y_1^iY_3^{j+h}Y_2^ju_\lambda \quad 0\leq j\leq m\quad 0\leq h\leq p\quad 0\leq i\leq j+n\quad 
j\neq 0\Rightarrow i\neq 0\right\}
$$
if  $\lambda=n\Lambda_1+m\Lambda_2+p\Lambda_3$ with $n>0$ and: \\
$
\mathfrak{G}_{\lambda}=\left\{Y_2^{m+h}Y_1Y_3^{1+h}Y_2u_\lambda,Y_2^{m+h}Y_3^{h}Y_2u_\lambda 
 \quad 0\leq h\leq p\right\}
$ 
if $\lambda=m\Lambda_2+p\Lambda_3$.\par\noindent 
 Although we do not need this fact in order to prove that  the $\mathfrak{a}_n$--module $V(\lambda)$ are indecomposable, let us first show that the set of generators $\mathfrak{G}$ is a minimal set of generators. We begin with 
\begin{lem}\label{g'g} No proper subset $\mathfrak{G}'_\lambda$ of $\mathfrak{G}_\lambda$ 
($\mathfrak{G}'_\lambda\subsetneq \mathfrak{G}_\lambda$) generates $\mathfrak{G}_\lambda$.\end{lem}
{\bf Proof}\  It suffices to show that any expression of the form
\beq\label{gg'eq}
\sum_{g\in \mathfrak{G}_\lambda}P_g(X_{2j},Y_{2j+1})Y^{(\mathbf{a_g})}u_\lambda\qquad P_g(X_{2j},Y_{2j+1})Y^{(\mathbf{a_g})}u_\lambda\neq 0\quad \forall g\in \mathfrak{G}_\lambda
\eeq
where  $P_{g}(X_{2j},Y_{2j+1})$ are non trivial polynomials  in the operators $X_{2j}$, $j=1,\cdots, 
\left[\frac{n}{2}\right]$, $Y_{2i+1}$, $i=1,\dots \left[\frac{n-1}{2}\right]$,
does not belong to the linear span $\langle \mathfrak{G_\lambda}\rangle$ of $\mathfrak{G_\lambda}$.\par\noindent
Let us denote 
by $V(\lambda)^-$ the linear span of all element of the Littelmann basis with $(a_k)^{2j}_{2j}<\lambda_{2j-1}^{2j}(H_{2j})$, since for every element of 
$\mathfrak{G}_\lambda$ yields $a^{2j}_{2j}=\lambda_{2j-1}^{2j}(H_{2j})$, we have \\  $V(\lambda)^-\cap \langle\mathfrak{G_\lambda}\rangle=\{0\}$. \par\noindent 
 Now from the proof of Theorem \ref{genset} point 1. for any element 
$Y^{(\mathbf{a})}u_\lambda$ in $\mathfrak{L}_\lambda$ and any operator $X_{2j}$, we have 
$X_{2j} Y^{(\mathbf{a})}u_\lambda=\sum_k c_k Y^{(\mathbf{a_k})}u_\lambda\in V(\lambda)^-$. Therefore it remains only to consider those combinations of the type  \rref{gg'eq} where there exists at least a monomial which contains only operators of odd index. For any such monomial  $P$  if 
 $V(\lambda)^+$ is  a subspace of  $V(\lambda)$  such that 
 $V(\lambda)= V(\lambda)^+\oplus (V(\lambda)^-\oplus\langle \mathfrak{G}_\lambda\rangle)$, then  from the proof of Theorem \ref{genset} points  3. and 4,  it follows that for any   $g\in \mathfrak{G}_\lambda$  $Pg=v^P_g+w^P_g$ with $v^P_g\in V(\lambda)^+$, $w^P_g\in (V(\lambda)^-\oplus\langle \mathfrak{G}_\lambda\rangle)$ and  $v^P_g\neq 0$, moreover if $g'\neq g$, $g,g'\in \mathfrak{G}_\lambda$,  or $P\neq Q$ then $v^P_g$ is linear independent from $v^Q_{g'}$. But then for any expression of type \rref{gg'eq} where there exist at least a monomial which is a product of only  the operators $Y_{2j+1}$  ($i=1,\dots \left[\frac{n-1}{2}\right]$) we have
$$
\sum_{g\in G}P_g(X_{2j},Y_{2j+1})Y^{(\mathbf{a_g})}u_\lambda\notin (V(\lambda)^-\oplus\langle \mathfrak{G}_\lambda\rangle)
$$
\endpf
\begin{theorem}\label{genminset} The set $\mathfrak{G}_\lambda$ is a minimal set of generators.
\end{theorem}
{\bf Proof} Let $\{w_1,\dots w_k\}$ be another set of generators, then for all 
$1\leq l\leq k$, choosing any ordering $G_\lambda=\{1,\dots,\#(\mathcal{G}_\lambda)\}$ (where $\#(S)$ denotes the number of elements in the set $S$) of the set $\mathcal{G}_\lambda$, we have:
$$
w_l=\sum_{g\in G_\lambda}a_{lg}Y^{(\mathbf{a_g})}u_\lambda+\sum_{g\in G_\lambda}P_{lg}(X_{2j},Y_{2j+1})Y^{(\mathbf{a_g})}u_\lambda
$$
where  $a_{lg}\in \mathbbm{C}$ and $P_{lg}
(X_{2j},Y_{2j+1})$ are polynomials  in the operators $X_{2j}$, $j=1,\cdots, 
\left[\frac{n}{2}\right]$, $Y_{2i+1}$, $i=1,\dots \left[\frac{n-1}{2}\right]$
  without constant term. Since the set  $\{w_1,\dots w_k\}$
generates $V(\Lambda)$  we may obtain acting on it the elements of 
$\mathfrak{G}_\lambda$. Let $\mathcal{T}=\{w_r \vert \ a_{rg}\neq 0 \mbox{ for some $g\in G_\lambda$}\}$ and
 $T=\{j\ \vert \ 1\leq j\leq k\ \vert\ w_r\in S\}$. Let $g\in G_\lambda$. Then   
$$
Y^{(\mathbf{a_g})}u_\lambda=\sum_{l\in T}b_{gl}\left(\sum_{g'\in G_\lambda}a_{lg'}Y^{(\mathbf{a_{g'}})}u_\lambda \right)
+\sum_{g'\in G_\lambda}P'_{tg'}(X_{2j},Y_{2j+1})Y^{(\mathbf{a_{g'}})}u_\lambda 
$$ 
where $P'_{tg'}$ are polynomials in the variables $X_{2j}$, $j=1,\cdots, 
\left[\frac{n}{2}\right]$, $Y_{2i+1}$, $i=1,\dots \left[\frac{n-1}{2}\right]$
 without constant term. From  Proof of Lemma \ref{g'g} it follows  that 
$\sum_{g'\in G_\lambda}P'_{tg'}(X_{2j},Y_{2j+1})Y^{(\mathbf{a_{g'}})}u_\lambda$ can not be equal to any combination 
of elements of $\mathfrak{G}_\lambda$. Hence 
$$
Y^{(\mathbf{a_g})}u_\lambda=\sum_{l\in T}b_{gl}\left(\sum_{g'\in G_\lambda}a_{lg'}Y^{(\mathbf{a_{g'}})}u_\lambda \right).
$$ 
This implies that if we put $B=\left(b_{gi}\right)_{\substack{g\in G_\lambda \\ i\in T}}$ and $A=\left(a_{ig}\right)_{\substack{i\in T \\ g\in G_\lambda}}$
 then $BA$ is the identity matrix. Hence $k\geq \#(T)\geq\mbox{rank}(B)\geq
 \#({G}_\lambda)$, so $\mathfrak{G}_\lambda$ is a minimal set of generators. The argument of this proof is due to Premat \cite{P}.
\endpf
\begin{cor}\label{genat} Let  $\mathfrak{W}=\{w_1,\dots w_k\}$ be a set  (non necessarily minimal) of generators, then there exist a injective map $\phi_\mathfrak{W}:\mathfrak{G}_\lambda\to \mathfrak{W}$, such that for every $Y^{(\mathbf{a_g})}u_\lambda\in \mathfrak{G}_\lambda$:
\beq\label{wgphi}
Y^{(\mathbf{a_g})}u_\lambda\mapsto w_{Y^{(\mathbf{a_g})}u_\lambda}=\phi_\mathfrak{W}(Y^{(\mathbf{a_g})}u_\lambda)=a_gY^{(\mathbf{a_g})}u_\lambda+\sum_{g'\in G_\lambda}P_{gg'}(X_{2j},Y_{2j+1})Y^{(\mathbf{a_{g'}})}u_\lambda
\eeq
for some $a_g\in \mathbbm{C}$, 
$a_g\neq 0$, where  $P_{gg'}$ are polynomials in the variables  $X_{2j}$, $j=1,\cdots, 
\left[\frac{n}{2}\right]$, $Y_{2i+1}$, $i=1,\dots \left[\frac{n-1}{2}\right]$, and the polynomial  $P_{gg}$ has no constant term.
\end{cor}
{\bf Proof} In the proof of Theorem \ref{genminset} we have shown that the elements of $\mathcal{W}$ can be written in the form 
$$
w_l=\sum_{g\in G_\lambda}a_{lg}Y^{(\mathbf{a_g})}u_\lambda
+\sum_{g'\in G_\lambda}P_{lg}(X_{2j},Y_{2j+1})Y^{(\mathbf{a_{g'}})}u_\lambda
$$
where $A=(a_{lt})_{\substack{l=1,\dots,\#(\mathfrak{W}) \\ t=1,\dots,\#({G}_\lambda)}}$ is a matrix of rank at least $\#(\mathcal{G}_\lambda)$ (recall that $\#(\mathfrak{W})\leq 
\#(\mathcal{G}_\lambda)$. This implies that for any $g\in 
\{1,\dots  \#(\mathcal{G}_\lambda)\}$ we can construct  a map $\phi:\{1,\dots  \#(\mathcal{G}_\lambda)\}\to \{1,\dots  \#(\mathcal{G}_\lambda)\}$ such that  for any $g\in 
\{1,\dots  \#(\mathcal{G}_\lambda)\}$, 
 $a_{\phi(g),g}$ is different from zero and $g\neq g'$ implies $\phi(g)\neq \phi(g')$. Then the map
$$
\begin{array}{ll}
\mathfrak{G}_\lambda&\to \mathfrak{W}\\
Y^{(\mathbf{a_g})}u_\lambda\mapsto&\phi_\mathfrak{W}(Y^{(\mathbf{a_g})}u_\lambda)=w_{\phi(g)}=
a_{\phi(g)g}Y^{(\mathbf{a_g})}u_\lambda+\sum_{\substack{g'\in G_\lambda \\  g\neq g'}}a_{\phi(g)g'}Y^{(\mathbf{a_{g'}})}u_\lambda\\&  +
\sum_{g'\in G_\lambda}P_{\phi(g)g'}(X_{2j},Y_{2j+1})(Y^{(\mathbf{a_{g'}})}u_\lambda)
\end{array}
$$
is  the searched  map. \endpf
\begin{prop}\label{g0} Any set  $\mathfrak{W}=\{w_1,\dots w_k\}$   of homogeneous generators contains an element $w_{{\overline{g}}}$ such that:
$$
w_{{\overline{g}}}=a_{{\overline{g}}}\overline{g}\qquad a_{{\overline{g}}}\neq 0 \in \mathbbm{C}.
$$
where  $\overline{g}$ is the element of the set of generators $\mathfrak{G}_\lambda$ given by
$$
\overline{g}=Y^{\left(\lambda^{2\left[\frac{n}{2}\right]}_{2\left[\frac{n}{2}\right]-1}
\left(H_{\left[\frac{n}{2}\right]}\right)\right)}_{2\left[\frac{n}{2}\right]}\cdots 
Y_{2j}^{(\lambda^{2j}_{2j-1}(H_{2j}))}\cdots  Y^{(\lambda^2_1(H_2))}_2u_\lambda.
$$
\end{prop}
{\bf Proof} 
The $\mathfrak{sl}(n+1,\mathbbm{C})$--weight $\mu_{\overline{g}}$ of the element $\overline{g}$ is  $\mu_{\overline{g}}=\lambda-\sum_{i=1}^{\left[\frac{n}{2}\right]}\lambda^{2i}_{2i-1}(H_{2i})\alpha_{2i}$ 
since 
$$
\mu_{\overline{g}}=s_{\alpha_{2\left[\frac{n}{2}\right]}}\cdots s_{\alpha_2}(\lambda)
$$ 
we have  $\mbox{dim}(V_{\mu_{\overline{g}}})=1$, but then the claim follows from Corollary \ref{genat} and the first part of this proposition.\endpf
Observe that we do not need the fact that $\mathfrak{G}_\lambda$ is a minimal set of generators in order to prove Proposition \ref{g0}. The simple fact that $\mathfrak{G}$ is a set of generators implies 
\begin{lem}\label{g0shway} Let  $\mathfrak{W}=\{w_1,\dots w_k\}$ be a set  (non necessarily minimal) of generators then there exists a $w_{\overline{k}}$ in  $\mathfrak{W}$ such that
$$
w_{\overline{k}}=\overline{a}\ \overline{g}+\sum_{g\in \mathfrak{G}_\lambda,g\neq \overline{g}}P_{lg}(X_{2j},Y_{2j+1})Y^{(\mathbf{a_g})}u_\lambda
$$ 
with $\overline{a}$ complex number different from zero.\end{lem}
{\bf Proof} Since:
$$
\overline{g}\notin X_{2j}(V(\lambda))\quad j=1,\cdots, 
\left[\frac{n}{2}\right] \qquad \overline{g}\notin Y_{2j+1}(V(\lambda))\quad 
i=1,\dots \left[\frac{n-1}{2}\right],
$$
the set $\mathfrak{W}$ is a  set of generators of $V(\lambda)$ only if it contains an element $w$ of the form
$$
w=\overline{a}\ \overline{g}+\sum_{g\in \mathfrak{G}_\lambda,g\neq \overline{g}}P_{lg}(X_{2j},Y_{2j+1})Y^{(\mathbf{a_g})}u_\lambda
$$ 
with $\overline{a}$ complex number different from zero.\endpf
Using Lemma \ref{g0shway} is obviously possible to prove directly Proposition \ref{g0}. \par
Let $\mathfrak{s}_n=\mathfrak{h}\rtimes \mathfrak{a}_n$ be  the subalgebra of 
$\mathfrak{sl}(n,\mathbbm{C})$ given by the semidirect product between the Cartan subalgebra $\mathfrak{h}$ and the subalgebra  $\mathfrak{a}_n$.  The $\mathfrak{sl}(n+1,\mathbbm{C})$--module ($\mathfrak{a}_n$--module) $V(\lambda)$ is also a  $\mathfrak{s}_n$--module, on which the subalgebra $\mathfrak{h}$ acts diagonally. Obviously  any set of generators of the $\mathfrak{a}_n$--module $V(\lambda)$  is also a set of generators of the $\mathfrak{s}_n$--module  $V(\lambda)$. Moreover for what said  above any  $\mathfrak{s}_n$--submodule of 
$V(\lambda)$ is a $\mathfrak{sl}(n+1,\mathbbm{C})$--weight module, i.e., it can decomposed as a direct sum of $\mathfrak{sl}(n+1,\mathbbm{C})$--weight spaces.   
From these facts it follows the 
\begin{prop}\label{enteneller} 
If the $\mathfrak{s}_n$--module $V(\lambda)$ decomposes in a direct sum of two subspaces: $V(\lambda)=U\oplus T$, then $\overline{g}$ belongs
 either to $U$ or to $T$.\end{prop}
{\bf Proof} Let  $\mathfrak{W}_U=\{w_1,\dots w_j\}$ and  $\mathfrak{W}_T=\{w_{j+1},\dots w_h\}$ be respectively a set of generators of $U$ and of $W$. Since  $U$ and $T$ are $\mathfrak{sl}(n+1,\mathbbm{C})$--weight modules we may suppose that both  $\mathfrak{W}_U$ and $\mathfrak{W}_T$ are made by homogeneous elements, and therefore 
$\mathfrak{W}=\mathfrak{W}_U\cup \mathfrak{W}_T=\{w_1,\dots w_j,w_{j+1},\dots w_h\}$  is a set of homogeneous generators of
$V(\lambda)$. Then  form proposition \ref{g0} it follows that there exists an index  $\overline{l}$,
$1\leq \overline{l}\leq h$ such that $\overline{g}=c_{\overline{g}}w_{\overline{l}}$. Hence $\overline{g}$ belongs either to $U$ or to $T$.\endpf 
\begin{theorem}\label{mthind} The $\mathfrak{a}_n$--module $V(\lambda)$ is indecomposable.\end{theorem}
{\bf Proof} Let us first show that  the  $\mathfrak{s}_n$--module $V(\lambda)$ is indecomposable.
Let us suppose that $V(\lambda)$ is the direct sum 
$V(\lambda)=U\oplus T$ of two   $\mathfrak{s}_n$--modules $U$ and $T$ and  let $\mathfrak{W}_U=\{w_1,\dots w_i\}$  (res.   $\mathfrak{W}_T=\{w_{i+1},\dots w_h\}$) be a set of homogeneous generators of $U$ (res. of $T$). We know from Proposition \ref{enteneller} that either $\overline{g}$ belongs to $U$ or to $T$. Say  $\overline{g}\in U$, then we shall show that $V(\lambda)=U$.\par
We say that an element  $Y^{\mathbf{(a)}}u_\lambda$ of the Littelmann basis  is of level $l$ if  $l$ is the minimal nonnegative integer such that $Y^{\mathbf{(a)}}u_\lambda=P_l\cdots P_1u_\lambda$ and any monomial $P_j$ $1\leq j\leq l$ is  a product of  elements $Y_i$ of index  either  odd or  even.\par 
It is immediate to see  that  all the elements of the Littelmann basis of length  $1$ and $0$ are in $\mathcal{U}(\mathfrak{a})(\overline{g})$ and therefore in $U$. Let now us suppose by induction that 
any element in $\mathfrak{G}_\lambda$  of level less or equal $l$  is  in $U$.  We need to  show that   
 any element in $\mathfrak{G}_\lambda$ of level $l+1$  also  belongs to $U$. First, since    any element $Y^{\mathbf{(a)}}u_\lambda$ in $\mathfrak{G}_\lambda$ is of the type $Y^{\mathbf{(a)}}u_\lambda=Y^{a^{2h}_{2h}}_{2h}\bigg(\cdots\bigg)u_\lambda$ 
with   $a^{2h}_{2h}\neq 0$, $0\leq h\leq \left[\frac{n}{2}\right]$,     
$\mathfrak{G}_\lambda$ decomposes as 
$$
\begin{array}{ll}
\mathfrak{G}_\lambda&=\bigcup_{1\leq j_1\leq \dots\leq j_s\leq \left[\frac{n}{2}\right]}\mathfrak{G}_{\lambda,j_1,\dots,j_s}\\
\mathfrak{G}_{\lambda,j_1,\dots,j_s}&=\left\{ g\in \mathfrak{G}_\lambda\ \vline  \begin{array}{ll} & g=Y^{(\mathbf{a})}u_\lambda=Y^{a^{2j_1}_{2j_1}}_{2j_1}(\cdots )Y^{a^{2j_s}_{2j_s}}_{2j_s}Y^{a^r_{2k+1}}_{2k_1}u_\lambda\\
&	\mbox{ with $a^{2j_i}_{2j_i}>0$  $i=1,\dots s$ ${a^r}_{2k+1}\neq 0$, $k<j_s$, $r>2j_s$}.\end{array}  \right\}
\end{array}
$$ 
Therefore it is enough to show that for any fixed set $\{j_1,\dots j_s\}$ ($1\leq j_1\leq \dots\leq j_s\leq \left[\frac{n}{2}\right]$)  the elements of length $l+1$ in $\mathfrak{G}_{\lambda,j_1,\dots,j_s}$  belong to $U$. We shall do it   by induction over  the orderings $\{\leq_{j_1},\dots,\leq_{j_s}\}$ defined in the proof of Theorem \ref{genset}.  If $g\in \mathfrak{G}_\lambda$ is minimal with respect all the ordering $\{\leq_{j_1},\dots,\leq_{j_s}\}$ then 
$$
X_{2j_1}^{a^{2j_1}_{2j_1}}\bigg(\cdots\bigg)X_{2j_s}^{a^{2j_s}_{2j_s}}g=c^{2j_1,\dots,2j_s}_{2j_1,\dots,2j_s}Y^{\mathbf{(a)}}u_\lambda=
c^{2j_1,\dots,2j_s}_{2j_1,\dots,2j_s}Y^{a^r_{2k+1}}_{2k+1}
\bigg(\cdots\bigg)u_\lambda
$$ 
with $c^{2j_1,\dots,2j_s}_{2j_1,\dots,2j_s}\neq 0$,  $a^r_{2k+1}\neq 0$ and $Y^{\mathbf{(a)}}u_\lambda=Y^{a^r_{2k+1}}_{2k+1}
\bigg(\cdots\bigg)u_\lambda$  element of  $\mathfrak{L}$ 
of level $l$. Since  $X_{2j_1}^{a^{2j_1}_{2j_1}}\bigg(\cdots\bigg)X_{2j_s}^{a^{2j_s}_{2j_s}}g$ has been obtained from an element of the set $\mathfrak{G}_\lambda$ of level $l+1$ by erasing the operators $Y_{2j_h}$, $1\leq h\leq s$,   and it is of the type
$Y^{a^r_{2k+1}}_{2k+1}
\bigg(\cdots\bigg)u_\lambda$ with $a^r_{2k+1}>0$, $0\leq k\leq \left[\frac{n-1}{2}\right]$, by the very definition of the set $\mathfrak{G}_\lambda$ it can be generated by an element of  $\mathfrak{G}_\lambda$ of level $l-1$ and it  is therefore by induction hypothesis a non trivial element in $U$. We can now  decompose $g$ as $g=g_U+g_T$ with $g_U=\sum_{k=1}^iP_k(X_{2s},Y_{2s+1})w_k\in U$ and $g_T
=\sum_{k=i+1}^hP_i(X_{2s},Y_{2s+1})w_k\in T$ and since all the elements $w_k$ are homogeneous, $g_U$ and $g_T$ are  of the same $\mathfrak{sl}(n+1,\mathbbm{C})$--weight of $g$. Now 
$$
X_{2j_1}^{a^{2j_1}_{2j_1}}\bigg(\cdots\bigg)X_{2j_s}^{a^{2j_s}_{2j_s}}g_T=
X_{2j_1}^{a^{2j_1}_{2j_1}}\bigg(\cdots\bigg)X_{2j_s}^{a^{2j_s}_{2j_s}}g-
X_{2j_1}^{a^{2j_1}_{2j_1}}\bigg(\cdots\bigg)X_{2j_s}^{a^{2j_s}_{2j_s}}g_U\in U\Longrightarrow X_{2j_1}^{a^{2j_1}_{2j_1}}\bigg(\cdots\bigg)X_{2j_s}^{a^{2j_s}_{2j_s}}g_T=0.
$$
 But the fact that $g_U$ and $g_T$ has the same  $\mathfrak{sl}(n+1,\mathbbm{C})$--weight of $g$ implies that 
they have also the same weight of $g$ with respect  any  subalgebra $\mathfrak{g}_{2j_r}$ spanned by the vector $H_{2j_r},X_{2j_r},Y_{2j_r}$ 
$1\leq r\leq s$ and equivalent to the complex simple Lie algebra $\mathfrak{sl}(2,\mathbbm{C})$.  Since 
$H_{2j_r}g=-a^{2j_r}_{2j_r}g$ with $a^{2j_r}_{2j_r}>0$  for $1\leq r\leq s$ the theory of the $\mathfrak{sl}(2,\mathbbm{C})$--finite dimensional modules implies that  for $1\leq r\leq s$,  $X_{2j_r}^{a^{2j_r}_{2j_r}}g_T=0$ if and only if $g_T=0$. Hence $g=g_U\in U$.  Now, since for any element $\widetilde{g}$ in  $\mathfrak{G}_{\lambda,j_1,\dots,j_s}$ which is not a minimal element for at least one of the ordering $\leq_{j_s}$ ($1\leq r\leq s$ ) we have  
$$
X_{2j_1}^{a^{2j_1}_{2j_1}}\bigg(\cdots\bigg)X_{2j_s}^{a^{2j_s}_{2j_s}}\widetilde{g}=\widetilde{c}^{2j_1,\dots,2j_s}_{2j_1,\dots,2j_s}
Y^{a^r_{2k+1}}_{2k+1}
\bigg(\cdots\bigg)u_\lambda+\sum_{\substack{(\mathbf{b})<_{j_r}(\mathbf{a})\\ s=1,\dots r}}\widetilde{c}_{\mathbf{(b)}}Y^{\mathbf{(b)}}u_\lambda
$$
by induction over the orderings $\leq_{j_s}$ ($1\leq r\leq s$ ) we have  that $X_{2j_1}^{a^{2j_1}_{2j_1}}\bigg(\cdots\bigg)X_{2j_s}^{a^{2j_s}_{2j_s}}\widetilde{g}\in U$. Then from the same  argument used above $\widetilde{g}\in U$. We have therefore proved that  any element of  $\mathfrak{G}_\lambda$ of length $l+1$ belong to $U$, if any  element of  $\mathfrak{G}_\lambda$ of length $l$ does. Therefore by induction the set of generators $\mathfrak{G}_\lambda$ belongs to $U$. Since $\mathfrak{G}_\lambda$ generates $V(\lambda)$ under the action of $\mathfrak{a}_n$, we have $V(\lambda)=U$ also as $\mathfrak{a}_n$--module. Hence the  $\mathfrak{a}_n$--module  $V(\lambda)$ is indecomposable for any integer dominant weight $\lambda$.\endpf


\end{document}